\documentclass[preprint,sort&compress,final,10pt]{elsarticle}
\usepackage{amsmath}
\usepackage{amssymb}
\usepackage{graphicx}
\usepackage{latexsym}
\usepackage{epstopdf}
\usepackage{natbib}
\usepackage{color}
\usepackage{hyperref}
\usepackage{amsmath,amscd}
 \usepackage[all]{xy}

\newtheorem{theorem}{Theorem}[section]

\newtheorem{lemma}[theorem]{Lemma}

\newtheorem{remark}[theorem]{Remark}
\numberwithin{equation}{section}

\def\Proof{\noindent{\bf Proof.}~}
\def\qed{\hfill$\square$\smallskip}
\def\dint{\displaystyle\int}
\def\dlim{\displaystyle\lim}

\def\dint{\displaystyle\int}

\def\dfrac#1#2{\frac{\displaystyle {#1}}{\displaystyle {#2}}}

\journal{\empty}
\date{}

\begin{document}

\begin{frontmatter}

\title{Boundedness in forced isochronous oscillators\footnote{Partially supported by the NSFC (11571041) and the Fundamental Research Funds for the Central Universities.}}

\author[au1]{Xiong Li\footnote{Corresponding author.}}

\address[au1]{School of Mathematical Sciences, Beijing Normal University, Beijing 100875, P.R. China.}

\ead[au1]{xli@bnu.edu.cn}

\author[au1]{Shasha Jin}

\ead[au1]{sjin@mail.bnu.edu.cn}

\begin{abstract}
In this paper we are concerned with the boundedness of all solutions for the forced isochronous oscillator
$$x''+V'(x)+g(x)=f(t),$$
where $V$ is a so-called $T$-isochronous potential, the perturbation $g$ is assumed to be bounded, and the $2\pi$-periodic function $f(t)$ is smooth. Using the resonant small twist theorem and averaged small twist theorem established by Ortega, we will prove the boundedness of all solutions for the above forced isochronous oscillator in the resonant and non-resonant cases under some reasonable assumptions, respectively.
\end{abstract}

\begin{keyword}
Boundedness; Forced isochronous oscillators; Resonant small twist theorem; Averaged small twist theorem.
\end{keyword}

\end{frontmatter}

\section{Introduction}

In this paper we are concerned with the boundedness of all solutions for the forced isochronous oscillator
\begin{equation}\label{FIO}
x''+V'(x)+g(x)=f(t),
\end{equation}
where $V$ is a so-called $T$-isochronous potential, the perturbation $g$ is bounded, and the $2\pi$-periodic function $f(t)$ is smooth. The origin  $(x,y)=(0,0)$ is called an isochronous center of the system
\begin{equation}\label{IO}
x'=y,\ \ \ \ y'=-V'(x),
\end{equation}
if every solution of system (\ref{IO}) is periodic with the minimal period $T>0$. Meanwhile, the equation
\begin{equation}\label{AIO}
x''+V'(x)=0
\end{equation}
is also called an isochronous system and $V$ is said to be a $T$-isochronous potential.

Obviously, the linear differential equation
$$x''+\omega^{2}x=0$$
is an isochronous system, since every solution is $2\pi/\omega$-periodic. In 1969, Lazer and Leach studied the existence of periodic solutions for the equation
$$x''+n^{2}x+g(x)=f(t)=f(t+2\pi), ~~n\in\mathbb{N}_+.$$
They showed that if the limits $\displaystyle\lim_{x\rightarrow\pm\infty}g(x):=g(\pm\infty)$ exist and are finite, and
$$\left|\int_{0}^{2\pi}f(t)e^{int}dt\right|<2\left|g(+\infty)-g(-\infty)\right|,$$
then this equation has at least one $2\pi$-periodic solution. Since then, the above inequality is called the Lazer-Landesman condition.

In 1999, Ortega \cite{Ortega99} studied a piecewise linear equation
\begin{equation}\label{PL}
x''+n^{2}x+h_{L}(x)=f(t), ~~n\in\mathbb{N}_+,
\end{equation}
where $f(t)\in C^{5}(\mathbb{R}/2\pi\mathbb{Z})$ and the piecewise linear function
$h_{L}(x)$ is given by
\begin{equation*}
h_{L}(x)=
\begin{array}{lll}
\left\{
\begin{array}{lll}
L,~~~~~x\geq1;\\[0.1cm]
Lx,~~~~|x|<1;\\[0.1cm]
-L,~~~~x\leq-1.\\
\end{array}
\right.
\end{array}
\end{equation*}
He proved that if
$$\left|\int_{0}^{2\pi}f(t)e^{int}dt\right|<4L,$$
then every solution of Eq. (\ref{PL}) is bounded, that is, if $x(t)$ is a solution of Eq. (\ref{PL}), then it exists on $(-\infty,+\infty)$ and
$$\sup_{t\in\mathbb{R}}(|x(t)|+|x'(t)|)<+\infty.$$

Liu\ \cite{Liu992} considered the general equation
\begin{equation}\label{GE}
x''+n^{2}x+\phi(x)=f(t),~~n\in\mathbb{N}_+,
\end{equation}
where $f(t)\in C^{7}(\mathbb{R}/2\pi\mathbb{Z})$, $\phi(x)\in C^{6}(\mathbb{R})$, the limits
$\lim \limits_{x\rightarrow\pm\infty}\phi(x)=\phi(\pm\infty)$
are finite and
$\lim \limits_{|x|\rightarrow+\infty}x^{6}\phi^{(6)}(x)=0.$
Then every solution of Eq. (\ref{GE}) is bounded if
$$\left|\int_{0}^{2\pi}f(t)e^{int}dt\right|<2|\phi(+\infty)-\phi(-\infty)|,$$
which is exactly the Lazer-Landesman condition. The above results demonstrate that the Lazer-Landesman condition also plays a key role in studying the boundedness problem.

The asymmetric oscillator
\begin{equation}\label{AO}
x''+ax^{+}-bx^{-}=0
\end{equation}
is also an isochronous system, where $x^{+}=\max(x,0)$, $x^{-}=\max(-x,0)$, $a, b$ are two different positive constants, since every solution of Eq. (\ref{AO}) is $2\pi/\omega$-periodic,
where $$\omega=2\left(\dfrac{1}{\sqrt{a}}+\dfrac{1}{\sqrt{b}}\right)^{-1}.$$
We remark that if $a=b=n^2$, then $\omega=n$. The forced asymmetric oscillator
\begin{equation}\label{a3}
x''+ax^{+}-bx^{-}=f(t)
\end{equation}
was first considered by Dancer \cite{Dancer76}, \cite{Dancer77} and Fu\v{c}ik \cite{Fucik80}.
They looked at this equation as a model of the so-called ``equations with jumping nonlinearities" and studied its periodic and Dirichlet boundary value problems. For recent developments, we refer to \cite{Gallouet82}, \cite{Habets93}, \cite{Lazer89}, \cite{Zhang98} and the references therein.

In 1996, Ortega\ \cite{Ortega96}  proved that all solutions of (\ref{a3}) are bounded if
$$f(t)=1+\varepsilon h(t),$$
where $h$ is smooth and $\varepsilon$ is small enough.  This result is in contrast with the
well-known phenomenon of linear resonance that occurs in the case $a=b=n^{2}$.

Liu\ \cite{Liu99} considered the boundedness of all solutions of Eq. (\ref{a3}) under the resonant case
$$\frac{1}{\sqrt{a}}+\frac{1}{\sqrt{b}}=\frac{2m}{n}\in \mathbb{Q}.$$
Let us recall this result. For a given $2\pi$-periodic function $f(t)$, define
\begin{equation}\label{Def}
\Phi_{f}(\theta)=\int_{0}^{2\pi}f(\theta+mt)C(mt)dt,
\end{equation}
and
$$\mathcal{A}(f)= \{\theta \in \mathbb{R}:\Phi_{f}(\theta) =0\},$$
where $C(t)$ is the solution of the equation
$$x''+ax^{+}-bx^{-}=0$$
with the initial conditions $x(0)=1, x'(0)=0$. He proved that if $\mathcal{A}(f)$ is empty, then all solutions of (\ref{a3}) are bounded.

On the other hand, Alonso and Ortega \cite{Alonso98} proved that if $\mathcal{A}(f)$ is not empty and
$$\Phi_{f}'(\theta)\not=0, ~~\mbox{for all}\quad \theta\in\mathcal{A}(f),$$
then all solutions of (\ref{a3}) with large initial conditions are unbounded. If
$$\frac{1}{\sqrt{a}}+\frac{1}{\sqrt{b}}\notin \mathbb{Q},$$
Ortega \cite{Ortega01} proved that if $f(t)\in C^{4}(\mathbb{R}/2\pi\mathbb{Z})$ and $[f]=\dfrac{1}{2\pi}\int_{0}^{2\pi}f(t)dt\neq0$, then all solutions of (\ref{a3}) are bounded.

In 2000, Fabry and Mawhin \cite{Fabry2000}, \cite{Fabry20} suggested to study the boundedness of all solutions for the equation
\begin{equation}\label{PAO}
x''+ax^{+}-bx^{-}+g(x)=f(t),
\end{equation}
where $a$ and $b$ are two positive constants, $g(x)$ is a bounded perturbation, and $f(t)$ is a smooth $2\pi$-periodic function. Wang \cite{Wang03} considered this question and obtained the boundedness of all solutions under some reasonable assumptions.

In 2009, Bonheure and Fabry \cite{Bonheure09} considered the boundedness of all solutions of the forced isochronous oscillator
$$
x''+V'(x)=f(t),
$$
where $V$ is a $T$-isochronous potential, $\lim \limits_{x\rightarrow+\infty}V''(x)=a>0,~~\lim \limits_{x\rightarrow-\infty}V''(x)=b>0$, $f$ is $T$-periodic, obtained the same result as that in \cite{Liu99}. Also they gave an example for such potential as
$$V(x)=\dfrac{x^{2}(x+2)^{2}}{4+2(1+\sigma)(x^{2}+2x)+4(x+1)\sqrt{1+\sigma x(x+2)}},$$
where $\sigma\in[0,1)$.

The above isochronous systems are defined on the whole real line. The following equation
$$x''+\frac{x+1}{4}-\frac{1}{4(x+1)^{3}}=0$$
is also an isochronous system, since all solutions are $2\pi$-periodic, and is not defined on $\mathbb{R}$, the potential tends to infinity as $x\rightarrow-1$. Liu \cite{Liu09}\ obtained the boundedness of all solutions of the forced isochronous oscillators with a repulsive singularity under the Lazer-Landesman condition. For more information and examples of isochronous centers, we refer to \cite{Chavarriga01} and the references therein.

Motivated by the above works, especially by \cite{Bonheure09} and \cite{Liu09}, in this paper we want to investigate the boundedness of all solutions for the forced isochronous oscillator (\ref{FIO}). Now we formulate our main result. Let $\omega: =2\pi/T$, where $T$ is the minimal period of solutions for the autonomous isochronous system (\ref{AIO}) , and $2\pi$ is the minimal period of the internal force $f(t)$.  We suppose that the following assumptions hold:\\

\noindent (1)\, $V(0)=V'(0)=0, V''(x)>0$ for $x\not =0$,
$$\lim \limits_{x\rightarrow0+}V''(x):=V''(0+)>0,~~~\lim \limits_{x\rightarrow0-}V''(x):=V''(0-)>0,$$
and
$$\lim \limits_{x\rightarrow+\infty}V''(x)=a>0,~~~\lim \limits_{x\rightarrow-\infty}V''(x)=b>0;$$\\
(2)\, $V(x)\in C^{6}(\mathbb{R}\setminus\{0\})$, the limits $\lim \limits_{x\rightarrow0+}V^{(k)}(x),\lim \limits_{x\rightarrow0-}V^{(k)}(x)$ exist and are finite for $3\leq k\leq 6$, and
$$\lim_{|x|\rightarrow+\infty}x^{4}V^{(6)}(x)=0;$$\\
(3)\, $g(x)\in C^{6}(\mathbb{R})$, the limits
$$\lim \limits_{x\rightarrow+\infty}g(x):=g(+\infty),~~\lim \limits_{x\rightarrow-\infty}g(x):=g(-\infty)$$
are finite and
$$\lim \limits_{|x|\rightarrow+\infty}x^{6}g^{(6)}(x)=0.$$

Then we have
\begin{theorem}\label{thm1.1}
Assume that $f(t)\in C^{6}(\mathbb{R}/2\pi\mathbb{Z})$ and the above hypotheses (1)-(3) hold. If\ $\omega\in\mathbb{Q}$, that is, there are two relatively prime positive integers  $m$, $n$ such that $\omega=\dfrac{n}{m}$ and
\begin{equation}\label{RC}
\sqrt{b}\cdot(\sqrt{a}+\sqrt{b})\cdot\Phi_{f}(\theta)\neq4\left[bg(+\infty)-ag(-\infty)\right],~~\theta\in\mathbb{R},
\end{equation}
then all solutions of Eq. (\ref{FIO}) are bounded; if\ $\omega\notin\mathbb{Q}$ and
\begin{equation}\label{NRC}
(b-a)[f]\neq[bg(+\infty)-ag(-\infty)],
\end{equation}
then all solutions of Eq. (\ref{FIO}) are bounded.
\end{theorem}

\begin{remark}
Firstly, from the hypothesis (1), there exist two positive constants $c_1, c_2$ such that $c_1 x^2\leq V(x)\leq c_2 x^2$, $c_1 x\leq V'(x)\leq c_2 x$ for all $x\in \mathbb{R}$. Also it follows from the hypothesis (2) that
$$\lim_{|x|\rightarrow+\infty}x^{k-2}V^{(k)}(x)=0, \ \ \ \ 3\leq k\leq 6.$$
Thus
$$\left|x^{k}V^{(k)}(x)\right|\leq C\,V(x),~~x\in\mathbb{R},~~1\leq k\leq 6,$$
where $C>0$ is a constant, $V^{(k)}(0)$ is understood as the limits $\lim \limits_{x\rightarrow0+}V^{(k)}(x)$, $\lim \limits_{x\rightarrow0-}V^{(k)}(x)$ for $k\geq 2$.

Define
\begin{equation}\label{a2}
W(x)£º=\frac{V(x)}{V'(x)},\ \  W(0):=\lim \limits_{x\rightarrow0}\frac{V(x)}{V'(x)}=\lim \limits_{x\rightarrow0}\frac{V'(x)}{V''(x)}=0,
\end{equation}
clearly $W(x)\in C^{6}(\mathbb{R}\setminus{0})$. Indeed, we have $W(x)\in C^{1}(\mathbb{R})$. For $x\not= 0$,
$$
W'(x)=1-\frac{V(x)V''(x)}{V'(x)^2},
$$
$$
\begin{array}{lll}
W'(0)&=&\lim \limits_{x\rightarrow0}\dfrac{V(x)}{xV'(x)}=\lim \limits_{x\rightarrow0}\dfrac{V'(x)}{V'(x)+xV''(x)}\\[0.4cm]
&=&\lim \limits_{x\rightarrow0}\dfrac{V''(x)}{2V''(x)+xV^{(3)}(x)}=\dfrac{1}{2},
\end{array}
$$
and
$$
\begin{array}{lll}
\lim \limits_{x\rightarrow0}W'(x)&=&1-\lim \limits_{x\rightarrow0}\dfrac{V(x)V''(x)}{V'(x)^2}
=1-\lim \limits_{x\rightarrow0}\dfrac{V'(x)V''(x)+V(x)V^{(3)}(x)}{2V'(x)V''(x)}\\[0.4cm]
&=&\dfrac{1}{2}=W'(0).
\end{array}
$$
Moreover, from the hypotheses (1) and (2), there also is a constant $C>0$ such that  for each $1\leq k\leq 6$,
\begin{equation}\label{a4}
|W(x)|\leq C|x|,~~|x^{k-1}W^{(k)}(x)|\leq C,~~ x\in\mathbb{R},
\end{equation}
where the value $x^{k-1}W^{(k)}(x)$ at $x=0$ is understood as the limit $\lim \limits_{x\rightarrow0}x^{k-1}W^{(k)}(x)$ for $k\geq 2$.

All the above estimates will be used to prove that $x(\theta,I)$ has the polynomial property, see Lemma \ref{lemma2.2} in Section 2. Similarly, it follows from the hypothesis (3) that for each \ $1\leq k\leq 6$,
\begin{equation}\label{a5}
\lim \limits_{|x|\rightarrow+\infty}x^{k}g^{(k)}(x)=0.
\end{equation}
\end{remark}
\begin{remark}
The proof of this theorem is based on the resonant small twist theorem (the resonant case: $\omega\in\mathbb{Q}$) and averaged small twist theorem (the non-resonant case: $\omega\not\in\mathbb{Q}$) established by Ortega \cite{Ortega99} and \cite{Ortega01}, respectively. The hypotheses (1)-(3) are used to prove that the Poincar\'{e} map of (\ref{FIO}) satisfies the assumptions of Ortega's theorems. Indeed, in the non-resonant case, we only need $f(t)\in C^{4}(\mathbb{R}/2\pi\mathbb{Z})$.
\end{remark}
\begin{remark}
When $V'(x)=ax^{+}-bx^{-}$, then Eq. (\ref{FIO}) takes the form (\ref{PAO}), which was investigated by Wang\, \cite{Wang03}. Although Eq. (\ref{FIO}) is more general than Eq. (\ref{PAO}), the results are completely same as that in \cite{Wang03}.  Since we can not introduce the explicit action and angle variables, we use some estimate methods similar to that in \cite{Liu09}.
\end{remark}
\begin{remark}
We would like to point out an interesting result of Ortega \cite{Ortega02}. In this paper, he showed
that there is a periodic function $p$ such that all solutions of the equation
$$x'' + V'(x) = \epsilon p(t)$$
are unbounded, where $V$ is an isochronous potential, $\epsilon$ is a small parameter. This result may show that the condition of Lazer-Landesman type (\ref{RC}) is necessary for the boundedness of all solutions.
\end{remark}

The rest of this paper is organized as follows. After introducing action and angle variables in Section 2, we state some technical lemmas, which will be used to prove our main result of the paper. Then we will give an asymptotic formula of the solutions of the autonomous isochronous system (\ref{AIO}). In Section 3, we will introduce another action and angle variables, and give an asymptotic expression of the Poincar\'{e} map. The main result will be proved by the resonant small twist theorem \cite{Ortega00} in Section 4 and averaged small twist theorem \cite{Ortega01} in Section 5, respectively.

\section{Action and angle variables}

In this section we first introduce action and angle variables. Let $y=x'$, then Eq. (\ref{FIO}) is equivalent to the following Hamiltonian system
\begin{equation}\label{b1}
x'=\frac{\partial H}{\partial y},~~~~~~y'=-\frac{\partial H}{\partial x},
\end{equation}
where the Hamiltonian is
$$H(x,y,t)=\frac{1}{2}y^{2}+V(x)+G(x)-xf(t)$$
with $G(x)=\int_{0}^{x}g(s)ds$.

In order to introduce action and angle variables, we consider the auxiliary autonomous system
$$x'=y,~~~~~~y'=-V'(x).$$
From our assumptions we know that all solutions of this system are $T$-periodic. For every $h>0$, denote by $I(h)$ the area enclosed by the closed curve
$$\dfrac{1}{2}y^{2}+V(x)=h.
$$
Let $-\alpha_{h}<0<\beta_{h}$ be such that $V(-\alpha_{h})=V(\beta_{h})=h$. Then by  hypotheses (1) it follows that
$$\lim \limits_{h\rightarrow+\infty}\alpha_{h}=\lim \limits_{h\rightarrow+\infty}\beta_{h}=+\infty.$$
In fact,
$$
\lim \limits_{h\rightarrow+\infty}\dfrac{\beta_h}{\sqrt{2a^{-1}h}}=1,~~~~\lim \limits_{h\rightarrow+\infty}\dfrac{\alpha_h}{\sqrt{2b^{-1}h}}=1.
$$
Moreover, it is easy to see that
$$I(h)=2\int_{-\alpha_{h}}^{\beta_{h}}\sqrt{2(h-V(s))}ds,~~ h>0.$$
Let
$$T(h)=2\int_{-\alpha_{h}}^{\beta_{h}}\frac{ds}{\sqrt{2(h-V(s))}},$$
$$T_{-}(h)=2\int_{-\alpha_{h}}^{0}\frac{ds}{\sqrt{2(h-V(s))}},$$
$$T_{+}(h)=2\int_{0}^{\beta_{h}}\frac{ds}{\sqrt{2(h-V(s))}},$$
then
$$I'(h)=T(h)=T_{-}(h)+T_{+}(h).$$
Since all solutions are $2\pi/\omega$-periodic, we have $T(h)=2\pi\omega^{-1}$, which yields that $I(h)=2\pi\omega^{-1}h$ and the inverse function of $I(h)$ is $h(I)=\dfrac{\omega }{2\pi}I$.

For every $(x,y)\in\mathbb{R}^{2}$, let us define the angle and action variables $(\theta,I)$ by
\begin{equation}\label{b2}
\theta(x,y)=
\begin{array}{ll}
\left\{
\begin{array}{ll}
\dint_{x}^{\beta_{h}}\frac{ds}{\sqrt{2(h(x,y)-V(s))}},~~~~~~~~~~~~~~~y\geq0,\\[0.6cm]
2\pi\omega^{-1}-\dint_{x}^{\beta_{h}}\frac{ds}{\sqrt{2(h(x,y)-V(s))}},~~~y<0,
\end{array}
\right.
\end{array}
\end{equation}
\begin{equation}\label{b3}
I(x,y)=2\dint_{-\alpha_{h}}^{\beta_{h}}\sqrt{2(h(x,y)-V(s))}ds,~~~~~~~~~~~~~~~~~~~~~~~~~
\end{equation}
where $h(x,y)=\frac{1}{2}y^{2}+V(x)$.

Obviously, the transformation $(\theta,I)\mapsto(x,y)$ is symplectic, thus (\ref{b1}) is transformed into another Hamiltonian system
\begin{equation}\label{b4}
\theta'=\frac{\partial H}{\partial I},~~~~I'=-\frac{\partial H}{\partial \theta},
\end{equation}
where the Hamiltonian
\begin{equation}\label{b5}
H(\theta,I,t)=I+2\pi\omega^{-1}G(x(\theta,I))-2\pi\omega^{-1} x(\theta,I)f(t)
\end{equation}
is $2\pi\omega^{-1}$ periodic with respect to $\theta$, $2\pi$ periodic with respect to $t$.

We first give the estimate on $x(\theta,I)$, whose proof is similar to that of Lemma A4.1 in \cite{Levi91}.
\begin{lemma}\label{lemma2.2}
There is a constant $C>0$ such that for $1\leq k\leq6$,
$$\left|I^{k}\partial_I^{k}x(\theta,I)\right|\leq C\,|x(\theta,I)|,$$
where $x=x(\theta,I)$ is defined implicitly by (\ref{b2}) and (\ref{b3}).
\end{lemma}
\Proof From the definition of $\theta$, we have, for $y\geq 0$,
$$
\begin{array}{lll}
\theta&=&\dint_{x}^{\beta_h}\frac{ds}{\sqrt{2(h-V(s))}}=\dint_{0}^{\beta_h}\frac{ds}{\sqrt{2(h-V(s))}}-\dint_{0}^{x}\frac{ds}{\sqrt{2(h-V(s))}}\\[0.8cm]
      &=&\dfrac{T_+(h)}{2}-\dint_{0}^{x}\frac{ds}{\sqrt{2(h-V(s))}}.
     \end{array}
$$
By the below Lemma \ref{lem4.1}, $\dfrac{T_+(h)}{2}=\dfrac{\pi}{2\sqrt{a}}$, taking the derivative with respect to the action variable $I$ in the both sides of the above equality
(the angle variable $\theta$ is independent of $I$) yields that
$$\partial_I\dint_{0}^{x}\frac{ds}{\sqrt{2(h-V(s))}}=0.$$
From \cite{Levi91} and \cite{Liu09}, one can get that
\begin{eqnarray*}
&&\partial_I\int_{0}^{x}\frac{ds}{\sqrt{2(h-V(s))}}\\
&=&\frac{1}{\sqrt{2(h-V(x))}}\left(\partial_I x-\frac{h'}{h}W(x)\right)+\frac{h'}{h}\int_{0}^{x}\left(\frac{1}{2}-W'(s)\right)\frac{ds}{\sqrt{2(h-V(s))}}.
\end{eqnarray*}
Since $h=h(I)=\dfrac{\omega }{2\pi}I$, thus
\begin{equation}\label{A3}
I\partial_I x=\sqrt{2(h-V(x))}\int_{0}^{x}\left(\frac{1}{2}-W'(s)\right)\frac{ds}{\sqrt{2(h-V(s))}}+W(x).
\end{equation}
For $0\leq s\leq x$, we have $\dfrac{h-V(x)}{h-V(s)}\leq 1 $, and by (\ref{a4}), we know that there exists $C>0$ such that
$$\left|I\partial_I x\right|\leq C\,|x|.$$

Using the properties on $V$ in Remark 1.2, the estimates for the derivatives of higher order and the case $y<0$ can be obtained in a same way as in \cite{Levi91} and we omit it here.\qed

Now we develop an asymptotic expression of $x(\theta,I)$ as $I\to+\infty$. First we define
\begin{equation*}
\Phi(x)=
\begin{array}{ll}
\left\{
\begin{array}{ll}
V(x)-\frac{a}{2}x^{2},~~~x\geq0,\\[0.2cm]
V(x)-\frac{b}{2}x^{2},~~~x<0.
\end{array}
\right.
\end{array}
\end{equation*}
By the assumptions (1) and (2), for $2\leq k\leq 6$, we have
\begin{equation}\label{d2}
\lim \limits_{|x|\rightarrow+\infty}x^{k-2}\Phi^{(k)}(x)=0.
\end{equation}

From the definition of $\theta$, it follows that
$$x_{\theta}(\theta,I)=-y(\theta,I).$$
Taking the derivative with respect to $\theta$ on both sides of the equation
$$\dfrac{1}{2}y^{2}+V(x)=h(I)=\dfrac{\omega }{2\pi}I$$ yields that
$$y_{\theta}(\theta,I)=V'(x),$$
which implies that
$$x_{\theta\theta}+V'(x)=0.$$

Define
$$\tilde{x}(\theta,I)=\frac{1}{\beta_{h}}x(\theta,I),$$
then
$$\tilde{x}(0,I)=1,~~~\tilde{x}_{\theta}(0,I)=0.$$
Obviously, if $\tilde{x}(\theta,I)\geq0$, $\tilde{x}(\theta,I)$ is the solution of the equation
\begin{equation}\label{d1}
\frac{d^{2}u}{d\theta^{2}}+au+\frac{1}{\beta_{h}}\Phi'(\beta_{h}u)=0
\end{equation}
with the initial conditions $u(0,I)=1, u_{\theta}(0,I)=0$; if $\tilde{x}(\theta,I)<0$, it is the solution of the equation
\begin{equation}\label{d11}
\frac{d^{2}u}{d\theta^{2}}+bu+\frac{1}{\beta_{h}}\Phi'(\beta_{h}u)=0
\end{equation}
with the initial conditions $u\left(\dfrac{T_{+}(h)}{2},I\right)=0, u_{\theta}\left(\dfrac{T_{+}(h)}{2},I\right)=-\dfrac{\sqrt{2h(I)}}{\beta_h}$.

By the definitions of $\theta$ and $\tilde{x}$, we also know that
$$\tilde{x}(\theta,I)>0~~\Leftrightarrow~~\theta\in\left(-\dfrac{T_{+}(h)}{2},\dfrac{T_{+}(h)}{2}\right),$$
$$\tilde{x}(\theta,I)=0~~\Leftrightarrow~~\theta=\dfrac{T_{+}(h)}{2},\ \  \dfrac{T_{+}(h)}{2}+T_{-}(h),$$
$$\tilde{x}(\theta,I)<0~~\Leftrightarrow~~\theta\in\left(\dfrac{T_{+}(h)}{2},\dfrac{T_{+}(h)}{2}+T_{-}(h)\right).$$
\begin{lemma}\label{lem4.1}
$\tilde{x}$\  has the following expression:
\begin{equation*}\label{d3}
\tilde{x}(\theta,I)=
\begin{array}{ll}
\left\{
  \begin{array}{ll}
    \cos\sqrt{a}\theta+X_{1}(\theta,I),~~~~~~~~~~~~~~~~~~~~~~~~~\theta\in\left[-\dfrac{\pi}{2\sqrt{a}},\dfrac{\pi}{2\sqrt{a}}\right],\\[0.4cm]
    -\sqrt{\dfrac{a}{b}}\sin\sqrt{b}\left(\theta-\dfrac{\pi}{2\sqrt{a}}\right)+X_{2}(\theta,I),~~~~\theta\in\left(\dfrac{\pi}{2\sqrt{a}},\dfrac{\pi}{2\sqrt{a}}+\dfrac{\pi}{\sqrt{b}}\right),
  \end{array}
\right.
\end{array}
\end{equation*}
where the functions $I^k\partial_I^k X_{1}(\theta,I)$ and $I^k\partial_I^k X_{2}(\theta,I)$ ($0\leq k\leq 6$) converge to $0$ uniformly for $\theta\in \left[-\dfrac{\pi}{2\sqrt{a}},\dfrac{\pi}{2\sqrt{a}}\right]$  and $\theta\in\left(\dfrac{\pi}{2\sqrt{a}},\dfrac{\pi}{2\sqrt{a}}+\dfrac{\pi}{\sqrt{b}}\right)$ as $I\rightarrow+\infty$, respectively.
\end{lemma}
\Proof When $\theta\in\left[-\frac{T_{+}(h)}{2},\frac{T_{+}(h)}{2}\right]$, \ $\tilde{x}$ is the solution of  (\ref{d1}) with the initial conditions $u(0,I)=1, ~u_{\theta}(0,I)=0$, thus
\begin{equation}\label{d111}
\tilde{x}(\theta,I)= \cos\sqrt{a}\theta-\frac{1}{\sqrt{a}}\int_{0}^{\theta}\frac{1}{\beta_{h}}\Phi'(\beta_{h}\tilde{x}(\tau,I))\sin\sqrt{a}(\theta-\tau)d\tau.
\end{equation}
Hence, the function $X_{1}$ is determined implicitly by
\begin{equation}\label{d1112}
X_{1}(\theta,I)=-\frac{1}{\sqrt{a}}\int_{0}^{\theta}\frac{1}{\beta_{h}}\Phi'(\beta_{h}(\cos\sqrt{a}\tau+X_{1}(\tau,I)))\sin\sqrt{a}(\theta-\tau)d\tau,
\end{equation}
where $\theta\in\left[-\frac{T_{+}(h)}{2},\frac{T_{+}(h)}{2}\right]$.

From (\ref{d2}), we know that
\begin{equation}\label{d0}
\lim\limits_{|x|\to+\infty}\dfrac{\Phi'(x)}{x}=0.
\end{equation}
Also, since
$$
\lim \limits_{I\rightarrow+\infty}\dfrac{\beta_h(I)}{\sqrt{\pi^{-1}a^{-1}\omega I}}=1,
$$
letting $I\rightarrow+\infty$ on both sides of (\ref{d1112}), by Lebesgue dominated theorem, the limit
$$\lim \limits_{I\rightarrow+\infty}X_{1}(\theta,I)=0$$
holds for any $\theta\in\left(-\frac{T_{+}(h)}{2},\frac{T_{+}(h)}{2}\right)$.

Now we are going to prove the above limit also holds uniformly for $\theta\in\left[-\frac{T_{+}(h)}{2},\frac{T_{+}(h)}{2}\right]$. Letting $I\to+\infty$ in (\ref{d111}) yields that
$$
\dlim_{I\to+\infty} \tilde{x}(\theta, I)=\cos\sqrt{a}\theta,\ \  \theta\in\left[-\frac{T_{+}(h)}{2},\frac{T_{+}(h)}{2}\right],
$$
also since $\tilde{x}\left(\frac{T_{+}(h)}{2}, I\right)=0$ for any $I$, therefore
$$
\cos\frac{\sqrt{a}\ T_{+}(h)}{2}=0,
$$
which implies that
$$
\dfrac{T_+(h)}{2}=\dfrac{\pi}{2\sqrt{a}}
$$
and
\begin{equation}\label{d00}
X_1\left(\dfrac{\pi}{2\sqrt{a}},I\right)=0.
\end{equation}

For any $\epsilon>0$, it follows from (\ref{d1112}) and (\ref{d0}) that $X_1(\theta, I)$ converges to $0$ uniformly for $\theta\in\left[-\frac{\pi}{2\sqrt{a}}+\epsilon,\frac{\pi}{2\sqrt{a}}-\epsilon\right]$ as $I\rightarrow+\infty$, which together with (\ref{d00}) and the continuity of $X_1(\theta, I)$ implies that the limit
$$\lim \limits_{I\rightarrow+\infty}X_{1}(\theta,I)=0$$
holds uniformly for $\theta\in\left[-\frac{\pi}{2\sqrt{a}},\frac{\pi}{2\sqrt{a}}\right]$.

Taking the derivative with respect to $I$ in the both sides of (\ref{d1112}), we can get that
$$
\partial_I X_{1}=-\frac{1}{\sqrt{a}}\int_{0}^{\theta}\left(-\frac{\Phi'}{\beta_{h}^{2}}\frac{d\beta_{h}}{dI}+\frac{\Phi''}{\beta_{h}}\frac{d\beta_{h}}{dI}(\cos\sqrt{a}\tau+X_{1})
+\Phi''\frac{\partial X_{1}}{\partial I}\right)\sin\sqrt{a}(\theta-\tau)d\tau.
$$
If we let
$$
a(\theta,I)=-\frac{I}{\sqrt{a}}\int_{0}^{\theta}\left(-\frac{\Phi'}{\beta_{h}^{2}}\frac{d\beta_{h}}{dI}
+\frac{\Phi''}{\beta_{h}}\frac{d\beta_{h}}{dI}(\cos\sqrt{a}\tau+X_{1})
\right)\sin\sqrt{a}(\theta-\tau)d\tau,
$$
and
$$
b(\theta,I,\tau)=-\frac{1}{\sqrt{a}}\Phi''\sin\sqrt{a}(\theta-\tau),
$$
then
$$
I\partial_I X_{1}=a(\theta,I)+\int_{0}^{\theta} b(\theta,I,\tau)I\partial_I X_{1}d\tau,
$$
and for $\theta>0$,
$$
\left|I\partial_I X_{1}\right|\leq |a(\theta,I)|+\int_{0}^{\theta} |b(\theta,I,\tau)|\left|I\partial_I X_{1}\right|d\tau.
$$
By Gronwall inequality, we have
$$\left|I\partial_I X_{1}\right|\leq |a(\theta,I)|+\int_{0}^{\theta} |a(\tau,I)||b(\theta,I,\tau)|\exp\left(\int_{\tau}^{\theta}|b(\theta,I,r)|dr\right)d\tau,$$
where $\theta>0$.

Since $|I\frac{d\beta_{h}}{dI}|\leq C \beta_{h}$, according to (\ref{d2}), (\ref{d1112}) and (\ref{d0}), for any $\epsilon>0$, $a(\theta, I)$ and $b(\theta,I,\tau)$ converges to $0$ uniformly for $\theta\in\left[0,\frac{\pi}{2\sqrt{a}}-\epsilon\right]$  and $\tau\in[0,\theta]$ as $I\rightarrow+\infty$, therefore  $I\partial_I X_{1}$ converges to $0$ uniformly for $\theta\in\left[0,\frac{\pi}{2\sqrt{a}}-\epsilon\right]$  as $I\rightarrow+\infty$.  Also from Lemma \ref{lemma2.2} we know that
$$
\partial_I X_{1}\left(\frac{\pi}{2\sqrt{a}},I\right)=0,
$$
which together with the continuity of $\partial_I X_{1}(\theta, I)$ implies that the limit
$$\lim \limits_{I\rightarrow+\infty}I\partial_I X_{1}(\theta, I)=0$$
holds uniformly for $\theta\in\left[0,\frac{\pi}{2\sqrt{a}}\right]$. According to the symmetry, the above limit also holds uniformly for $\theta\in\left[-\frac{\pi}{2\sqrt{a}},\frac{\pi}{2\sqrt{a}}\right]$. Differentiating (\ref{d1112}) with respect to $I$ repeatedly, the estimates for the derivatives of higher order can be obtained in a similar way.

When $\dfrac{\pi}{2\sqrt{a}}<\theta<\dfrac{\pi}{2\sqrt{a}}+T_{-}(h)$, then $\tilde{x}(\theta,I)<0$, and
\begin{equation}\label{d01}
\tilde{x}\left(\dfrac{\pi}{2\sqrt{a}}+T_{-}(h),I\right)=0,
\mathcal{}\end{equation}
and it is the solution of (\ref{d11}) with the initial conditions $u\left(\dfrac{\pi}{2\sqrt{a}},I\right)=0, u_{\theta}\left(\dfrac{\pi}{2\sqrt{a}},I\right)=-\dfrac{\sqrt{2h(I)}}{\beta_h}.$ Therefore,
\begin{equation}\label{dd}
\begin{array}{ll}
\tilde{x}(\theta,I)=&-\dfrac{\sqrt{2b^{-1}h(I)}}{\beta_h}\sin\sqrt{b}\left(\theta-\dfrac{\pi}{2\sqrt{a}}\right)\\[0.3cm]
&-\dfrac{1}{\sqrt{b}}\dint_{\dfrac{\pi}{2\sqrt{a}}}^{\theta}
\frac{1}{\beta_{h}}\Phi'(\beta_{h}\tilde{x}(\tau,I))\sin\sqrt{b}\left(\theta-\tau\right)d\tau,
\end{array}
\end{equation}
where $\theta\in\left[-\dfrac{\pi}{2\sqrt{a}},\dfrac{\pi}{2\sqrt{a}}+T_{-}(h)\right]$.

Since
$$
\lim \limits_{I\rightarrow+\infty}\dfrac{\sqrt{2b^{-1}h(I)}}{\beta_h}=\sqrt{\dfrac{a}{b}},
$$
by Lebesgue dominated theorem, letting $I\to+\infty$ in (\ref{dd}), we know that
$$
\dlim_{I\to+\infty} \tilde{x}(\theta, I)=-\sqrt{\dfrac{a}{b}}\sin\sqrt{b}\left(\theta-\dfrac{\pi}{2\sqrt{a}}\right),\ \  \theta\in\left[-\dfrac{\pi}{2\sqrt{a}},\dfrac{\pi}{2\sqrt{a}}+T_{-}(h)\right],
$$
which together with (\ref{d01}) implies that
$$
T_{-}(h)=\dfrac{\pi}{\sqrt{b}}, \ \ \tilde{x}\left(\dfrac{\pi}{2\sqrt{a}}+\dfrac{\pi}{\sqrt{b}},I\right)=0,
$$
and
$$
T(h)=\dfrac{\pi}{\sqrt{a}}+\dfrac{\pi}{\sqrt{b}},~~~~~\omega=2\left(\dfrac{1}{\sqrt{a}}+\dfrac{1}{\sqrt{b}}\right)^{-1}.
$$
Thus, we rewrite $\tilde{x}(\theta,I)$ as
\begin{equation}\label{d1111}
\begin{array}{ll}
\tilde{x}(\theta,I)=&-\sqrt{\dfrac{a}{b}}\sin\sqrt{b}\left(\theta-\dfrac{\pi}{2\sqrt{a}}\right)\\[0.4cm]
&+\left[\sqrt{\dfrac{a}{b}}-\dfrac{\sqrt{2b^{-1}h(I)}}{\beta_h}\right]\sin\sqrt{b}\left(\theta-\dfrac{\pi}{2\sqrt{a}}\right)\\[0.4cm]
&-\dfrac{1}{\sqrt{b}}\dint_{\dfrac{\pi}{2\sqrt{a}}}^{\theta}
\frac{1}{\beta_{h}}\Phi'(\beta_{h}\tilde{x}(\tau,I))\sin\sqrt{b}\left(\theta-\tau\right)d\tau,
\end{array}
\end{equation}
where $\theta\in\left[\dfrac{\pi}{2\sqrt{a}},\dfrac{\pi}{2\sqrt{a}}+\dfrac{\pi}{\sqrt{b}}\right]$, and the function $X_{2}$ is determined implicitly by
\begin{equation*}\label{d1111}
\begin{array}{ll}
X_{2}(\theta,I)=&\left[\sqrt{\dfrac{a}{b}}-\dfrac{\sqrt{2b^{-1}h(I)}}{\beta_h}\right]\sin\sqrt{b}\left(\theta-\dfrac{\pi}{2\sqrt{a}}\right)\\[0.4cm]
&-\dfrac{1}{\sqrt{b}}\dint_{\dfrac{\pi}{2\sqrt{a}}}^{\theta}
\frac{1}{\beta_{h}}\Phi'\left(\beta_{h}\left(-\sqrt{\dfrac{a}{b}}\sin\sqrt{b}\left(\tau-\dfrac{\pi}{2\sqrt{a}}\right)+X_2(\tau,I)\right)\right)\sin\sqrt{b}\left(\theta-\tau\right)d\tau.
\end{array}
\end{equation*}
Similar to the estimate on $X_1$, $I^k\partial_I^k X_{2}(\theta,I)$ ($0\leq k\leq 6$) converges to $0$ uniformly for $\theta\in\left[\dfrac{\pi}{2\sqrt{a}},\dfrac{\pi}{2\sqrt{a}}+\dfrac{\pi}{\sqrt{b}}\right]$ as $I\rightarrow+\infty$. Thus we have finished the proof of the lemma.\qed

Then we have
\begin{equation*}\label{d3}
x(\theta,I)=
\begin{array}{ll}
\left\{
  \begin{array}{ll}
    \beta_h\cos\sqrt{a}\theta+\beta_h X_{1}(\theta,I),~~~~~~~~~~~~~~~~~~~~~~~~~\theta\in\left[-\dfrac{\pi}{2\sqrt{a}},\dfrac{\pi}{2\sqrt{a}}\right],\\[0.4cm]
    -\beta_h\sqrt{\dfrac{a}{b}}\sin\sqrt{b}\left(\theta-\dfrac{\pi}{2\sqrt{a}}\right)+\beta_h X_{2}(\theta,I),~~~~\theta\in\left(\dfrac{\pi}{2\sqrt{a}},\dfrac{\pi}{2\sqrt{a}}+\dfrac{\pi}{\sqrt{b}}\right).
  \end{array}
\right.
\end{array}
\end{equation*}
Since
\begin{equation}\label{Beta}
\lim \limits_{I\rightarrow+\infty}\dfrac{\beta_h(I)}{\sqrt{\pi^{-1}a^{-1}\omega I}}=1,
\end{equation}
then
\begin{equation*}\label{d3}
x(\theta,I)=
\begin{array}{ll}
\left\{
  \begin{array}{ll}
   \sqrt{\pi^{-1}a^{-1}\omega}\,I^{\frac{1}{2}}\cos\sqrt{a}\theta+ \widetilde{X}_{1}(\theta,I),~~~~~~~~~~~~~~~~~~~~\theta\in\left[-\dfrac{\pi}{2\sqrt{a}},\dfrac{\pi}{2\sqrt{a}}\right],\\[0.4cm]
    - \sqrt{\pi^{-1}b^{-1}\omega}\,I^{\frac{1}{2}}\sin\sqrt{b}\left(\theta-\dfrac{\pi}{2\sqrt{a}}\right)+ \widetilde{X}_{2}(\theta,I),~~~~\theta\in\left(\dfrac{\pi}{2\sqrt{a}},\dfrac{\pi}{2\sqrt{a}}+\dfrac{\pi}{\sqrt{b}}\right),
  \end{array}
\right.
\end{array}
\end{equation*}
where the functions $\widetilde{X}_{1}$ and $\widetilde{X}_{2}$ are given by
\begin{equation}\label{X1}
\widetilde{X}_{1}=\left(\beta_h-\sqrt{\pi^{-1}a^{-1}\omega}\,I^{\frac{1}{2}}\right)\cos\sqrt{a}\theta+\beta_h X_1(\theta,I),
\end{equation}
\begin{equation}\label{X2}
\widetilde{X}_{2}=\left(\sqrt{\pi^{-1}b^{-1}\omega}\,I^{\frac{1}{2}}-\beta_h\sqrt{\dfrac{a}{b}}\right)\sin\sqrt{b}\left(\theta-\dfrac{\pi}{2\sqrt{a}}\right)+\beta_h X_2(\theta,I).
\end{equation}

For the sake of convenience, we denote the approximate expression of $x(\theta,I)$ by
\begin{equation*}\label{d3}
\bar{x}(\theta,I)=
\begin{array}{ll}
\left\{
  \begin{array}{ll}
   \sqrt{\pi^{-1}a^{-1}\omega}\,I^{\frac{1}{2}}\cos\sqrt{a}\theta,~~~~~~~~~~~~~~~~~~~~~\theta\in\left[-\dfrac{\pi}{2\sqrt{a}},\dfrac{\pi}{2\sqrt{a}}\right],\\[0.4cm]
    -\sqrt{\pi^{-1}b^{-1}\omega}\,I^{\frac{1}{2}}\sin\sqrt{b}\left(\theta-\dfrac{\pi}{2\sqrt{a}}\right),~~~~\theta\in\left(\dfrac{\pi}{2\sqrt{a}},\dfrac{\pi}{2\sqrt{a}}+\dfrac{\pi}{\sqrt{b}}\right).
  \end{array}
\right.
\end{array}
\end{equation*}
Moreover, if we assume that $C(\theta)$ is the solution of
$$x''+ax^{+}-bx^{-}=0$$
with the initial conditions $x(0)=1, x'(0)=0$, that is,
\begin{equation*}
C(\theta)=
\begin{array}{ll}
\left\{
  \begin{array}{ll}
   \cos\sqrt{a}\theta,~~~~~~~~~~~~~~~~~~~~~~~~~~~\theta\in\left[-\dfrac{\pi}{2\sqrt{a}},\dfrac{\pi}{2\sqrt{a}}\right],\\[0.4cm]
   -\sqrt{\dfrac{a}{b}}\sin\sqrt{b}\left(\theta-\dfrac{\pi}{2\sqrt{a}}\right),~~~~\theta\in\left(\dfrac{\pi}{2\sqrt{a}},\dfrac{\pi}{2\sqrt{a}}+\dfrac{\pi}{\sqrt{b}}\right),
  \end{array}
\right.
\end{array}
\end{equation*}
then
\begin{equation}\label{E1}
\bar{x}(\theta,I)=\sqrt{\pi^{-1}a^{-1}\omega}\,I^{\frac{1}{2}} C(\theta),
\end{equation}
and
\begin{equation}\label{E2}
x(\theta,I)=\bar{x}(\theta,I)+X(\theta,I),
\end{equation}
where
\begin{equation*}\label{d3}
X(\theta,I)=
\begin{array}{ll}
\left\{
  \begin{array}{ll}
   \widetilde{X}_{1}(\theta,I),~~~~\theta\in\left[-\dfrac{\pi}{2\sqrt{a}},\dfrac{\pi}{2\sqrt{a}}\right],\\[0.4cm]
   \widetilde{X}_{2}(\theta,I),~~~~\theta\in\left(\dfrac{\pi}{2\sqrt{a}},\dfrac{\pi}{2\sqrt{a}}+\dfrac{\pi}{\sqrt{b}}\right),
  \end{array}
\right.
\end{array}
\end{equation*}
and the limits
\begin{equation}\label{X2}
\lim \limits_{I\rightarrow+\infty}I^{k-\frac{1}{2}}\partial_I^{k} X=0, \ \ \ \ k=0,1,\cdots, 6
\end{equation}
hold uniformly for $\theta\in \left[-\dfrac{\pi}{2\sqrt{a}},\dfrac{\pi}{2\sqrt{a}}+\dfrac{\pi}{\sqrt{b}}\right]$.

\section{Another action and angle variables}

In this section we introduce another canonical transformation such that the transformed system is
a small perturbation of an integrable system. Now we go back to system (\ref{b4}). Observe that
$$Id\theta-Hdt=-(Hdt-Id\theta),$$
this means that if one can solve $I=I(t, H,\theta)$ from (\ref{b4}) as a function of $H$ ($\theta$ and $t$ as parameters), then
\begin{equation}\label{e1}
\frac{dH}{d\theta}=-\partial_t I(t,H,\theta),~~\frac{dt}{d\theta}=\partial_H I(t,H,\theta).
\end{equation}
That is, (\ref{e1}) is a Hamiltonian system with the Hamilton $I=I(t,H,\theta)$ and now the new action, angle and time variables are $H$, $t$ and $\theta$, respectively. The relation between (\ref{b4}) and (\ref{e1}) is that if $(I(t),\theta(t))$ is a solution of (\ref{b4}) and the inverse function $t(\theta)$ of $\theta(t)$ exists, then $(H(\theta,I(t(\theta)),t(\theta)),t(\theta))$ is a solution of (\ref{e1}) and vice versa.

Recall that
$$
H(\theta,I,t)=I+2\pi\omega^{-1}G(x(\theta,I))-2\pi\omega^{-1}x(\theta,I)f(t).
$$
Let
$$\Psi(\theta,I,t):=2\pi\omega^{-1}G(x(\theta,I))-2\pi\omega^{-1}x(\theta,I)f(t),$$
then
$$
H(\theta,I,t)=I+\Psi(\theta,I,t),
$$
and by the assumption (3) and Lemma \ref{lemma2.2}, there is a constant $C>0$ such
that for $k+l\leq6$,
$$\left|I^{k}\partial_I^k\partial_t^l\Psi(\theta,I,t)\right|\leq C\sqrt{I}.$$
Thus
$$\partial_I H=1+\partial_I\Psi\rightarrow 1, ~~~~I\rightarrow+\infty.$$
Hence, by the implicit function theorem, there is a function $R(t,H,\theta)$ such that
$$I=H-R(t,H,\theta),$$
where
$$R(t,H,\theta)=\Psi(\theta,H-R,t).$$
It is easy to see that
$$\left|H^{k}\partial_H^k \partial_t^lR(t, H,\theta)\right|\leq C\sqrt{H},$$
where $C>0$ is a constant. Furthermore,  if we let $R_{1}(t,H,\theta):=\Psi(\theta,H,t)-R(t,H,\theta)$, then
$$R_{1}(t,H,\theta)=-\int_{0}^{1}\partial_I\Psi(\theta,H-sR,t)\,Rds,$$
and there exists a positive constant $C$ such that for $k+l\leq6$,
\begin{equation}\label{e4}
\left|H^{k}\partial_H^k \partial_t^lR_{1}(t,H,\theta)\right|\leq C.
\end{equation}

The new Hamilton is written in the form
\begin{equation*}
\begin{array}{lll}
I&=& H- \Psi(\theta,H,t)+R_{1}(t,H,\theta)\\[0.3cm]
&=& H-2\pi\omega^{-1} G(x(\theta,H))+2\pi\omega^{-1} x(\theta,H)f(t)+R_{1}(t,H,\theta)
\end{array}
\end{equation*}
and system (\ref{e1}) is
\begin{equation*}
\begin{array}{ll}
\left\{
  \begin{array}{ll}
   \dfrac{dt}{d\theta}=\partial_H\, I=1-2\pi\omega^{-1}\partial_H\, x(\theta,H)\left[g(x (\theta,H))-f(t)\right]+\partial_H\, R_{1}(t,H,\theta), \\[0.4cm]
   \dfrac{dH}{d\theta}=-\partial_t\, I=-2\pi\omega^{-1} x(\theta,H)f'(t)-\partial_t\, R_{1}(t,H,\theta).
  \end{array}
\right.
\end{array}
\end{equation*}
Now we replace $\theta$ by $\omega^{-1}\theta$, then the system becomes
\begin{equation}\label{e3}
\left\{
  \begin{array}{ll}
   \dfrac{dt}{d\theta}=&\omega^{-1}-2\pi\omega^{-2}\partial_H\, x(\omega^{-1}\theta,H)\left[g(x(\omega^{-1}\theta,H))-f(t)\right]\\[0.2cm]
   &+\omega^{-1}\partial_H\, R_{1}(t,H,\omega^{-1}\theta),\\[0.2cm]
   \dfrac{dH}{d\theta}=&-2\pi\omega^{-2} x(\omega^{-1}\theta,H)f'(t)-\omega^{-1}\partial_t\, R_{1}(t,H,\omega^{-1}\theta),
  \end{array}
\right.
\end{equation}
which is $2\pi$ periodic with respect to $t$ and $\theta$, respectively.

Introduce a new action variable $\rho\in[1,2]$ and a parameter $\epsilon>0$ by $H=\epsilon^{-2}\rho$. Then, $H\gg1\Leftrightarrow0<\epsilon\ll1$. Under this transformation, system (\ref{e3}) is changed into the form
\begin{equation}\label{e5}
\left\{
  \begin{array}{ll}
   \dfrac{dt}{d\theta}=&\omega^{-1}-2\pi\omega^{-2}\partial_H\, x(\omega^{-1}\theta,\epsilon^{-2}\rho)\left[g(x(\omega^{-1}\theta,\epsilon^{-2}\rho))-f(t)\right]\\[0.2cm]&+\omega^{-1}\partial_H\, R_{1}(t,\epsilon^{-2}\rho,\omega^{-1}\theta),\\[0.2cm]
   \dfrac{d\rho}{d\theta}=&-2\pi\omega^{-2}\epsilon^{2} x(\omega^{-1}\theta,\epsilon^{-2}\rho)f'(t)-\omega^{-1}\epsilon^{2}\partial_t\, R_{1}(t,\epsilon^{-2}\rho,\omega^{-1}\theta),
  \end{array}
\right.
\end{equation}
which is also the Hamiltonian system with the Hamilton
\begin{eqnarray*}
\begin{array}{lll}
\Gamma(t,\rho,\theta;\epsilon)
&=&\omega^{-1}\rho-2\pi\omega^{-2}\epsilon^{2}[G(x(\omega^{-1}\theta,\epsilon^{-2}\rho))-x(\omega^{-1}\theta,\epsilon^{-2}\rho)f(t)]\\[0.2cm]
&&+\omega^{-1}\epsilon^{2} R_{1}(t,\epsilon^{-2}\rho,\omega^{-1}\theta).
\end{array}
\end{eqnarray*}

Obviously, if $0<\epsilon\ll1$, the solution $(t(\theta,t_{0},\rho_{0}),\rho(\theta,t_{0},\rho_{0}))$ of (\ref{e5}) with the initial data $(t_{0},\rho_{0})\in\mathbb{R}\times[1,2]$ is defined in the interval $\theta\in [0,2\pi]$ and $\rho(\theta,t_{0},\rho_{0})\in\left[1/2,3\right]$ for $\theta\in [0,2\pi]$. Hence the Poincar\'{e} map of  (\ref{e5}) is well defined in the domain $\mathbb{R}\times[1,2]$, and has the intersection property (see \cite{Ortega00}).

From now on, we use the notations $o_{k}(1)$ and $O_{k}(1)$. A function $f(t,\rho,\theta;\epsilon)$ is said to be of order $o_{k}(1)$ if it is $C^{k}$ in $(t,\rho)$ and for $k_{1}+k_{2}\leq k$,
$$\lim \limits_{\epsilon\rightarrow0}\left|\partial_t^{k_{1}}\partial_\rho^{k_{2}}f(t,\rho,\theta;\epsilon)\right|=0,~~\mbox{uniformly~in}~(t,\rho,\theta).$$
We say a function $f(t,\rho,\theta;\epsilon)\in O_{k}(1)$ if $f(t,\rho,\theta;\epsilon)\in C^{k}$ in $(t,\rho)$ and for $k_{1}+k_{2}\leq k$,
$$\left|\partial_t^{k_{1}}\partial_\rho^{k_{2}}f(t,\rho,\theta;\epsilon)\right|\leq C,$$
where $C>0$ is a constant independent of the arguments $t,\rho,\theta,\epsilon$.

Now we first give some estimates, which will be used to calculate the asymptotic expression of the Poincar\'{e} map of (\ref{e5}) as $\epsilon\ll1$.
Suppose that the solution of (\ref{e5}) with the initial condition $(t(0),\rho(0))=(t_{0},\rho_{0})$ is of the form
$$t=t_{0}+\omega^{-1}\theta+\epsilon\,\Sigma_{1}(t_{0},\rho_{0},\theta;\epsilon),~~~~\rho=\rho_{0}+\epsilon\,\Sigma_{2}(t_{0},\rho_{0},\theta;\epsilon).$$
Then the Poincar\'{e} map $P$ of (\ref{e5}) is
$$P: ~~t_{1}=t_{0}+2\pi\omega^{-1}+\epsilon\,\Sigma_{1}(t_{0},\rho_{0},2\pi;\epsilon),~~~~\rho_{1}=\rho_{0}+\epsilon\,\Sigma_{2}(t_{0},\rho_{0},2\pi;\epsilon),$$
and the functions $\Sigma_1$ and $\Sigma_2$ satisfy
$$
\begin{array}{lll}
\Sigma_{1}&=&-2\pi\omega^{-2}\epsilon^{-1}\dint_{0}^{\theta}\partial_H x(\omega^{-1}\theta,\epsilon^{-2}\rho)[g(x(\omega^{-1}\theta,\epsilon^{-2}\rho))-f(t)]d\theta\\[0.3cm]
&&+\omega^{-1}\epsilon^{-1}\dint_{0}^{\theta}\partial_H R_1(t,\epsilon^{-2}\rho,\omega^{-1}\theta)d\theta,\\[0.3cm]
\Sigma_{2}&=&-2\pi\omega^{-2}\epsilon\dint_{0}^{\theta}x(\omega^{-1}\theta,\epsilon^{-2}\rho)f'(t)d\theta-\omega^{-1}\epsilon
\dint_{0}^{\theta}\partial_t R_1(t,\epsilon^{-2}\rho,\omega^{-1}\theta)d\theta,
\end{array}
$$
where $t=t_{0}+\omega^{-1}\theta+\epsilon\,\Sigma_{1}, \rho=\rho_{0}+\epsilon\,\Sigma_{2}.$

By Lemma \ref{lemma2.2}, (\ref{e4}) and the assumptions (1)-(3), we know that the terms in the right-hand side of the above equations are bounded, that is,
$$|\Sigma_{1}|+|\Sigma_{2}|\leq C, ~~~\theta\in[0,2\pi],$$
where $C>0$ is a constant. Hence, for $\rho_{0}\in[1,2]$, we may choose $\epsilon$ sufficiently small such that
\begin{equation}\label{f1}
\rho_{0}+\epsilon\Sigma_{2}\geq\frac{\rho_{0}}{2}\geq\frac{1}{2},~~~(t_{0},\theta)\in[0,2\pi]\times[0,2\pi].
\end{equation}
 Similar to the proof in \cite{Dieckerhoff87}, one can obtain
\begin{equation}\label{f2}
\Sigma_{1}\in O_{6}(1),~~~~\Sigma_{2}\in O_{5}(1).
\end{equation}

\begin{lemma}\label{lem6.1}
The following estimates hold:
$$x(\omega^{-1}\theta,\epsilon^{-2}\rho)-x(\omega^{-1}\theta,\epsilon^{-2}\rho_{0})\in O_{6}(1),~~~~~~~~~~~$$
$$\partial_H x(\omega^{-1}\theta, \epsilon^{-2}\rho)-\partial_H x(\omega^{-1}\theta,\epsilon^{-2}\rho_{0})\in \epsilon^{2}O_{5}(1).$$
\end{lemma}
\Proof Let
\begin{equation}\label{Delta}
\Delta(t_0,\rho_0,\theta;\epsilon):=x(\omega^{-1}\theta,\epsilon^{-2}\rho)-x(\omega^{-1}\theta,\epsilon^{-2}\rho_{0})
=\dint_0^1 \partial_H x(\omega^{-1}\theta,\epsilon^{-2}\rho_{0}+s\epsilon^{-1}\Sigma_2)\epsilon^{-1}\Sigma_2ds.
\end{equation}
By Lemma \ref{lemma2.2}, (\ref{f1}), (\ref{f2}), we have 
$$
|\Delta(t_0,\rho_0,\theta;\epsilon)|\leq C\dfrac{\epsilon^{-1}\Sigma_2}{\sqrt{\epsilon^{-2}\rho_{0}+s\epsilon^{-1}\Sigma_2}}\leq C.
$$
Take the derivative with respect to $\rho_0$ in the both sides of (\ref{Delta}), we have
$$
\partial_{\rho_0} \Delta=\dint_0^1\partial_H^2 x\cdot(\epsilon^{-2}+s\epsilon^{-1}\partial_{\rho_0}\Sigma_2)\epsilon^{-1}\Sigma_2ds.
$$
Using Lemma \ref{lemma2.2}, (\ref{f2}), one may find a constant $C>9$ such that $|\partial_{\rho_0} \Delta|\leq C$. Analogously, one may obtain, by a direct but cumbersome computation, that
$$
|\partial_{t_0}^{k_1}\partial_{\rho_0}^{k_2}\Delta(t_0,\rho_0,\theta;\epsilon)|\leq C
$$
for $k_1+k_2\leq 6$. The estimates for $\partial_H x(\omega^{-1}\theta, \epsilon^{-2}\rho)-\partial_H x(\omega^{-1}\theta,\epsilon^{-2}\rho_{0})$ 
follow from a similar argument, we omit it here. \qed

\begin{lemma}\label{lem6.2}
The following estimate holds:
$$
\epsilon^{-1}\int_{0}^{2\pi}\partial_H\bar{x}(\omega^{-1}\theta,\epsilon^{-2} \rho_{0})g(x(\omega^{-1}\theta,\epsilon^{-2}\rho_{0}))d\theta
=\frac{\omega^{\frac{3}{2}}}{\sqrt{\pi\rho_{0}}}\left(\frac{1}{a}g(+\infty)-\frac{1}{b}g(-\infty)\right)+o_{5}(1).
$$
\end{lemma}
\Proof Let
$$
\bar{g}(\rho_0;\epsilon):=\epsilon^{-1}\int_{0}^{2\pi}\partial_H\bar{x}(\omega^{-1}\theta,\epsilon^{-2} \rho_{0})g(x(\omega^{-1}\theta,\epsilon^{-2}\rho_{0}))d\theta.
$$ 
Recall that
$$
\overline{x}(\theta,I)=\sqrt{\pi^{-1}a^{-1}\omega}\,I^{\frac{1}{2}} C(\theta),
$$
and
\begin{equation*}
C(\theta)=
\left\{
  \begin{array}{ll}
   \cos\sqrt{a}\theta,~~~~~~~~~~~~~~~~~~~~~~~~~~~\theta\in\left[-\dfrac{\pi}{2\sqrt{a}},\dfrac{\pi}{2\sqrt{a}}\right],\\[0.4cm]
   -\sqrt{\dfrac{a}{b}}\sin\sqrt{b}\left(\theta-\dfrac{\pi}{2\sqrt{a}}\right),~~~~\theta\in\left(\dfrac{\pi}{2\sqrt{a}},\dfrac{\pi}{2\sqrt{a}}+\dfrac{\pi}{\sqrt{b}}\right).
  \end{array}
\right.
\end{equation*}
Therefore, we obtain
\begin{eqnarray*}
\bar{g}(\rho_0;\epsilon)
&=&\epsilon^{-1}\,\omega\int_{0}^{\frac{2\pi}{\omega}}\partial_H\bar{x}(\theta,\epsilon^{-2} \rho_{0})g(x(\theta,\epsilon^{-2}\rho_{0}))d\theta\\
&=&\frac{\omega^{\frac{3}{2}}}{2\sqrt{a\pi\rho_{0}}}\int_{0}^{\frac{2\pi}{\omega}}g(x(\theta,\epsilon^{-2}\rho_{0}))C(\theta)d\theta.
\end{eqnarray*}
By Lebesgue dominated theorem, we have
\begin{eqnarray*}
\lim \limits_{\epsilon\rightarrow0^+}\bar{g}(\rho_0;\epsilon)
&=&\frac{\omega^{\frac{3}{2}}}{\sqrt{a\pi\rho_{0}}}\int_{0}^{\frac{\pi}{2\sqrt{a}}}g(+\infty)\cos\sqrt{a}\theta d\theta\\
&&-\frac{\omega^{\frac{3}{2}}}{\sqrt{a\pi\rho_{0}}}\sqrt{\dfrac{a}{b}}\int_{\frac{\pi}{2\sqrt{a}}}^{\frac{\pi}{\omega}}g(-\infty)\sin\sqrt{b}\left(\theta-\frac{\pi}{2\sqrt{a}}\right)d\theta\\
&=&\frac{\omega^{\frac{3}{2}}}{\sqrt{\pi\rho_{0}}}\left(\frac{1}{a}g(+\infty)-\frac{1}{b}g(-\infty)\right).
\end{eqnarray*}

Since
$$
\partial_{\rho_0}\bar{g}(\rho_0;\epsilon)=-\dfrac{1}{2\rho_0}\bar{g}(\rho_0;\epsilon)
+\frac{\omega^{\frac{3}{2}}}{2\sqrt{a\pi\rho_{0}}}\int_{0}^{\frac{2\pi}{\omega}}g'(x(\theta,\epsilon^{-2}\rho_{0}))\partial_H x(\theta,\epsilon^{-2}\rho_{0})\epsilon^{-2} C(\theta)d\theta,
$$
by the assumption (\ref{a5}), Lemma \ref{lemma2.2} and Lebesgue dominated theorem, we know that 
$$
\lim \limits_{\epsilon\rightarrow0^+}\int_{0}^{\frac{2\pi}{\omega}}g'(x(\theta,\epsilon^{-2}\rho_{0}))\partial_H x(\theta,\epsilon^{-2}\rho_{0})\epsilon^{-2} C(\theta)d\theta=0,
$$
and 
$$
\lim \limits_{\epsilon\rightarrow0^+}\partial_{\rho_0}\bar{g}(\rho_0;\epsilon)
=-\frac{\omega^{\frac{3}{2}}}{2\sqrt{\pi\rho_{0}^3}}\left(\frac{1}{a}g(+\infty)-\frac{1}{b}g(-\infty)\right).
$$
The estimates for the derivatives of higher order can be obtained in a similar way.\qed

\section{The resonant case}
In this section we will prove the main result under the resonant case:\ $\omega\in\mathbb{Q}$, that is, there are two relatively prime positive integers $m$, $n$ such that $\omega=\dfrac{n}{m}$. Introducing the new time variable by $\theta=n\vartheta$, then the corresponding Hamiltonian system is
\begin{equation}\label{f4}
\left\{
  \begin{array}{ll}
   \dfrac{dt}{d\vartheta}=m-2m\pi\omega^{-1}\partial_H\, x(m\vartheta,\epsilon^{-2}\rho)\left[g(x(m\vartheta,\epsilon^{-2}\rho))-f(t)\right]+m\partial_H R_{1}, \\[0.4cm]
   \dfrac{d\rho}{d\vartheta}=-2m\pi\omega^{-1}\epsilon^{2} x(m\vartheta,\epsilon^{-2}\rho)f'(t)-m\epsilon^{2}\partial_t R_{1},
  \end{array}
\right.
\end{equation}
where $R_{1}=R_{1}(t,\epsilon^{-2}\rho,m\vartheta)$.

We assume that the  solution of (\ref{f4}) with the initial condition $(t(0),\rho(0))=(t_{0},\rho_{0})$ is of
the form
$$t=t_{0}+m\vartheta+\epsilon f_{1}(t_{0},\rho_{0},\vartheta;\epsilon),~~~~\rho=\rho_{0}+\epsilon f_{2}(t_{0},\rho_{0},\vartheta;\epsilon),$$
where the functions $f_1$ and $f_2$ satisfy
$$
\begin{array}{lll}
f_{1}=&-2m\pi\omega^{-1}\epsilon^{-1}\dint_{0}^{\vartheta}\partial_H\, x(m\vartheta,\epsilon^{-2}\rho)\left[g(x(m\vartheta,\epsilon^{-2}\rho))-f(t)\right]d\vartheta\\[0.3cm]
&+m\epsilon^{-1}\dint_{0}^{\vartheta}\partial_H R_{1}d\vartheta,\\[0.3cm]
f_{2}=&-2m\pi\omega^{-1}\epsilon\dint_{0}^{\vartheta}x(m\vartheta,\epsilon^{-2}\rho)f'(t)d\vartheta-m\epsilon\dint_{0}^{\vartheta}\partial_t R_{1}d\vartheta,
\end{array}
$$
and $t=t_{0}+m\vartheta+\epsilon f_{1},\rho=\rho_{0}+\epsilon f_{2}.$
Then, the Poincar\'{e} map of (\ref{f4}) is
$$P:~~t_{1}=t_{0}+2m\pi+\epsilon f_{1}(t_{0},\rho_{0},2\pi;\epsilon),~~~~\rho=\rho_{0}+\epsilon f_{2}(t_{0},\rho_{0},2\pi;\epsilon).$$

By (\ref{a5}), (\ref{X2}), (\ref{e4}), (\ref{f2}), and  Lemmas \ref{lem6.1},  \ref{lem6.2}, we can get
\begin{eqnarray*}
&&f_{1}(t_{0},\rho_{0},2\pi;\epsilon)\\
&=&-2m\pi\omega^{-1}\epsilon^{-1}\int_{0}^{2\pi}\partial_H x(m\vartheta,\epsilon^{-2}(\rho_{0}+\epsilon f_{2}))g(x(m\vartheta,\epsilon^{-2}(\rho_{0}+\epsilon f_{2})))d\vartheta\\
&&+2m\pi\omega^{-1}\epsilon^{-1}\int_{0}^{2\pi}\partial_H x(m\vartheta,\epsilon^{-2}(\rho_{0}+\epsilon f_{2}))f(t_{0}+m\vartheta+\epsilon f_{1})d\vartheta+\epsilon O_{5}(1)\\
&=&-2m\pi\omega^{-1}\epsilon^{-1}\int_{0}^{2\pi}\partial_H x(m\vartheta,\epsilon^{-2}\rho_{0})g(x(m\vartheta,\epsilon^{-2}\rho_{0}))d\vartheta\\
&&+2m\pi\omega^{-1}\epsilon^{-1}\int_{0}^{2\pi}\partial_H x(m\vartheta,\epsilon^{-2}\rho_{0})f(t_{0}+m\vartheta)d\vartheta+\epsilon O_{5}(1)\\
&=&-2m\pi\omega^{-1}\epsilon^{-1}\int_{0}^{2\pi}\partial_H \bar{x}(m\vartheta,\epsilon^{-2}\rho_{0})g(x(m\vartheta,\epsilon^{-2}\rho_{0}))d\vartheta\\
&&+2m\pi\omega^{-1}\epsilon^{-1}\int_{0}^{2\pi}\partial_H \bar{x}(m\vartheta,\epsilon^{-2}\rho_{0})f(t_{0}+m\vartheta)d\vartheta+\epsilon O_{5}(1)\\
&=&-2m\sqrt{\frac{\pi\omega}{\rho_0}}\left(\frac{g(+\infty)}{a}-\frac{g(-\infty)}{b}\right)
+m\sqrt{\frac{\pi}{a\omega\rho_0}}\int_{0}^{2\pi}f(t_{0}+m\vartheta)C(m\vartheta)d\vartheta+\epsilon O_{5}(1)
\end{eqnarray*}
and
\begin{eqnarray*}
&&f_{2}(t_{0},\rho_{0},2\pi;\epsilon)\\
&=&-2m\pi\omega^{-1}\epsilon\int_{0}^{2\pi}x(m\vartheta,\epsilon^{-2}(\rho_{0}+\epsilon f_{2}))f'(t_{0}+m\vartheta+\epsilon f_{1})d\vartheta+\epsilon O_{5}(1)\\
&=&-2m\pi\omega^{-1}\epsilon\int_{0}^{2\pi}x(m\vartheta,\epsilon^{-2}\rho_{0})f'(t_{0}+m\vartheta)d\vartheta+\epsilon O_{5}(1)\\
&=&-2m\pi\omega^{-1}\epsilon\int_{0}^{2\pi}\bar{x}(m\vartheta,\epsilon^{-2}\rho_{0})f'(t_{0}+m\vartheta)d\vartheta+\epsilon O_{5}(1)\\
&=&-2m\sqrt{\frac{\pi\rho_0}{a\omega}}\int_{0}^{2\pi} C(m\vartheta)f'(t_{0}+m\vartheta)d\vartheta+\epsilon O_{5}(1).
\end{eqnarray*}

Hence the Poincar\'{e} map has the form
\begin{equation}\label{f5}
\begin{array}{ll}
P:\left\{
\begin{array}{ll}
t_{1}=t_{0}+2m\pi-\epsilon m \pi^{\frac{1}{2}}\omega^{-\frac{1}{2}}l_{1}(t_{0})\rho_0^{-\frac{1}{2}}+\epsilon o_{5}(1), \\[0.4cm]
\rho_{1}=\rho_{0}-2\epsilon m\pi^{\frac{1}{2}}a^{-\frac{1}{2}}\omega^{-\frac{1}{2}} l_{2}(t_{0})\rho_{0}^{\frac{1}{2}}+\epsilon o_{5}(1),
\end{array}
\right.
\end{array}
\end{equation}
where
$$
\begin{array}{ll}
l_{1}(t_{0})&=2\omega\left(\dfrac{g(+\infty)}{a}-\dfrac{g(-\infty)}{b}\right)
-\dfrac{1}{\sqrt{a}}\dint_{0}^{2\pi}f(t_{0}+m\vartheta)C(m\vartheta)d\vartheta,\\[0.4cm]
l_{2}(t_{0})&=\dint_{0}^{2\pi}f'(t_{0}+m\vartheta)C(m\vartheta)d\vartheta.
\end{array}
$$
Under the diffeomorphism
$$t=t,~~ r=\frac{1}{\rho},$$
the map $P$ is transformed into the following form
\begin{equation}\label{f6}
\begin{array}{ll}
\overline{P}:\left\{
\begin{array}{ll}
t_{1}=t_{0}+2m\pi-\epsilon m \pi^{\frac{1}{2}}\omega^{-\frac{1}{2}}l_{1}(t_{0})r_{0}^{\frac{1}{2}}+\epsilon o_{5}(1), \\[0.4cm]
r_{1}=r_{0}+2\epsilon m\pi^{\frac{1}{2}}a^{-\frac{1}{2}}\omega^{-\frac{1}{2}} l_{2}(t_{0})r_{0}^{\frac{3}{2}}+\epsilon o_{5}(1).
\end{array}
\right.
\end{array}
\end{equation}
If $l_{1}(t_{0})\neq0$, that is, for any $t_0\in \mathbb{R}$,
$$2\omega\left(\dfrac{1}{a}g(+\infty)-\dfrac{1}{b}g(-\infty)\right)\neq\dfrac{1}{\sqrt{a}}\int_{0}^{2\pi}f(t_{0}+m\vartheta)C(m\vartheta)d\vartheta,$$
same as in \cite{Wang03}, it is easy to verify that (\ref{f6}) satisfied all assumptions of the resonant small twist theorem in \cite{Ortega99}. Thus, all solutions of (\ref{FIO}) are bounded.

\section{The non-resonant case}
In this section we will prove the main result under the non-resonant case:\ $\omega\notin\mathbb{Q}$. Similar to the resonant case, one can obtain that the expression of the Poincar\'{e} map is
\begin{equation}\label{f7}
\begin{array}{ll}
P:\left\{
\begin{array}{ll}
t_{1}=t_{0}+2\pi\omega^{-1}-\epsilon\Sigma_{1}(t_{0},\rho_{0},2\pi;\epsilon)+\epsilon o_{4}(1), \\[0.4cm]
\rho_{1}=\rho_{0}-\epsilon\Sigma_{2}(t_{0},\rho_{0},2\pi;\epsilon)+\epsilon o_{4}(1),
\end{array}
\right.
\end{array}
\end{equation}
where
$$\Sigma_{1}(t_{0},\rho_{0},2\pi;\epsilon)=\omega^{-\frac{3}{2}}\sqrt{\dfrac{\pi}{\rho_0}}\left[\dfrac{2\omega}{a}g(+\infty)-\dfrac{2\omega}{b}g(-\infty)
-\dfrac{1}{\sqrt{a}}\int_{0}^{2\pi}C(\omega^{-1}\theta)f(t_{0}+\omega^{-1}\theta)d\theta\right],$$
$$\Sigma_{2}(t_{0},\rho_{0},2\pi;\epsilon)=
2\omega^{-\frac{3}{2}}\sqrt{\dfrac{\pi\rho_{0}}{a}}\int_{0}^{2\pi}C(\omega^{-1}\theta)f'(t_{0}+\omega^{-1}\theta)d\theta.~~~~~~~~~~~~~~~~~~~~~~~~~~~~~~~~~~~$$
Thus we have
\begin{eqnarray*}
&&\int_{0}^{2\pi}\dfrac{\partial\Sigma_{1}}{\partial\rho_{0}}(t_{0},\rho_{0},2\pi;\epsilon)dt_{0}\\
&=&-\int_{0}^{2\pi}\frac{\sqrt{\pi}}{2}\omega^{-\frac{3}{2}}\rho_{0}^{-\frac{3}{2}}\left[\dfrac{2\omega}{a}g(+\infty)-\dfrac{2\omega}{b}g(-\infty)
-\dfrac{1}{\sqrt{a}}\int_{0}^{2\pi}C(\omega^{-1}\theta)f(t_{0}+\omega^{-1}\theta)d\theta\right]dt_{0}\\
&=&-\frac{\sqrt{\pi}\omega^{-\frac{3}{2}}\rho_{0}^{-\frac{3}{2}}}{2}\left[4\pi\omega\left(\dfrac{g(+\infty)}{a}-\dfrac{g(-\infty)}{b}\right)
-\dfrac{1}{\sqrt{a}}\int_{0}^{2\pi}C(\omega^{-1}\theta)\int_{0}^{2\pi}f(t_{0}+\omega^{-1}\theta)dt_{0}d\theta\right]\\
&=&-\frac{\sqrt{\pi}\omega^{-\frac{3}{2}}\rho_{0}^{-\frac{3}{2}}}{2}\left[4\pi\omega\left(\dfrac{g(+\infty)}{a}-\dfrac{g(-\infty)}{b}\right)
-\dfrac{2\pi}{\sqrt{a}}[f]\int_{0}^{2\pi}C(\omega^{-1}\theta)d\theta\right]\\
&=&-2\pi^{\frac{3}{2}}\omega^{-\frac{1}{2}}\rho_{0}^{-\frac{3}{2}}\left[\left(\dfrac{g(+\infty)}{a}-\dfrac{g(-\infty)}{b}\right)-[f]\left(\frac{1}{a}-\frac{1}{b}\right)\right],
\end{eqnarray*}
where $[f]=\dfrac{1}{2\pi}\dint_{0}^{2\pi}f(t)dt$.
If
$$bg(+\infty)-ag(-\infty)\neq[f](b-a),$$
same as in \cite{Wang03}, it is easy to verify that (\ref{f7}) satisfied all assumptions of the averaged small twist theorem in \cite{Ortega01}. Therefore, all solutions of (\ref{FIO}) are bounded.

\section*{References}
\bibliographystyle{elsarticle-num}

\end{document}